\magnification=\magstep1
\input amssym
\centerline{\bf ON A THEOREM OF KN\"ORR}
\bigskip\noindent
\centerline{Burkhard K\"ulshammer, Institute for Mathematics,}
\centerline{Friedrich Schiller University, Jena, Germany}
\bigskip\bigskip\noindent
{\narrower\narrower {\bf Abstract.} Kn\"orr has constructed an ideal, 
in the center of the 
$p$-modular group algebra of a finite group $G$, whose dimension is
the number of $p$-blocks of defect zero in $G/Q$; here $p$ is a 
prime and $Q$ is a normal $p$-subgroup of $G$. We generalize his
construction to symmetric algebras.\par}
\bigskip\noindent
Let $F$ be an algebraically closed field of characteristic $p>0$, 
and let $Q$ be a normal $p$-subgroup of a finite group $G$. Motivated
by Alperin's Weight Conjecture, Kn\"orr in [3] proved that the number
of blocks of defect zero in the group algebra $F[G/Q]$ can be computed
as the dimension of a certain ideal in the center ${\rm Z}(FG)$. 

Kn\"orr's approach used the permutation $FG$-module on the cosets of a
Sylow $p$-subgroup $P$ of $G$ and its endomorphism ring. Here we 
present a different approach to his result by making use of properties
of symmetric algebras. Applied to group algebras, our main result is
as follows:
\bigskip\noindent
{\bf Theorem 1.} {\it In the situation described above, the number of
blocks of defect zero in $F[G/Q]$ coincides with the dimension of the 
ideal ${\rm W}_Q(FG) := {\rm Tr}_Q^G(FG \cdot G_p^+)$ of ${\rm Z}(FG)$.}
\bigskip\noindent 
Let us explain the notation used in Theorem 1. We denote by $G_p$ the 
set of $p$-elements in $G$, and by $G_{p'}$ the set of $p'$-elements
in $G$. For any subset $X$ of $G$, we set $X^+ := \sum_{x \in X} x \in
FG$. Then $FG \cdot G_p^+ = G_p^+ \cdot FG$ is a principal (two-sided)
ideal of $FG$ contained in the (left and right) socle ${\rm S}(FG)$. 
By a result of Tsushima in [8] (see also [4]), the annihilator of 
$G_p^+$ in $FG$ is the sum of the radical ${\rm J}(FG)$ of $FG$ and the 
left ideals $FGe$ of $FG$ where $e$ ranges over the primitive idempotents
in $FG$ such that $\dim FGe$ is divisible by $p|P|$. 

For any subgroup $X$ of $G$,
$$(FG)^X := \{a \in FG: xax^{-1} = a \hbox{ for all } x \in X\}$$
is a subalgebra of $FG$ and, for any subgroup $Y$ of $X$,
$${\rm Tr}_Y^X: (FG)^Y \longrightarrow (FG)^X, \quad a \longmapsto
\sum_{xY \in X/Y} xax^{-1},$$
is the relative trace (transfer) map. Observe that, in Theorem 1, we
have $FG \cdot G_p^+ \subseteq (FG)^Q$ since $uG_p^+ = G_p^+$ for all
$u \in Q$. Since $FG \cdot G_p^+$ is an ideal in $(FG)^Q$, standard
properties of the transfer map imply that ${\rm W}_Q(FG) := 
{\rm Tr}_Q^G(FG \cdot G_p^+)$ is an ideal in $(FG)^G = {\rm Z}(FG)$. 

By a result of Robinson [7], the number of blocks of defect zero of a
group algebra can be computed as the rank of a certain matrix with
coefficients in ${\Bbb Z}/p{\Bbb Z}$; see also [5]. In [9], Wang and
Zhang have translated Robinson's theorem from $G/Q$ to $G$. 

This paper is organized as follows: In Section 1, we first recall some
properties of symmetric algebras and then prove a version of Kn\"orr's
result in this context. In Section 2 we apply the results of Section 1
to group algebras, prove Theorem 1 and finish with some additional 
remarks.
\vfill\eject
\bigskip\bigskip\noindent
{\bf 1. Symmetric algebras}
\bigskip\noindent
Let $F$ be an algebraically closed field, and let $A$ be a symmetric 
$F$-algebra with symmetrizing linear form $\lambda: A \longrightarrow 
F$. Moreover, let $I$ be an ideal of $A$ such that the $F$-algebra 
$A/I$ is also symmetric, and let $\mu: A/I \longrightarrow F$ be a 
corresponding symmetrizing linear form. Then there exists a unique 
element $z$ in the center ${\rm Z}(A)$ of $A$ such that $\mu(a+I)
= \lambda(az)$ for all $a \in A$. Consequently, $I$ is the annihilator
of $z$ in $A$. We denote by $\nu: A \longrightarrow A/I$, $a \longmapsto
a+I$, the canonical epimorphism, and by $\nu^\ast: A/I \longrightarrow
A$ the adjoint of $\nu$ (with respect to $\lambda$ and $\mu$). Thus
$$\lambda(\nu^\ast(x)a) = \mu(x \nu(a)) \hbox{ for all } x \in A/I, \;
a \in A.$$
Explicitly, we have
$$\nu^\ast(a+I) = az \hbox{ for all } a \in A;$$
in particular, $\nu^\ast$ is a monomorphism of $A$-$A$-bimodules. This
implies that $\nu^\ast({\rm J}(A/I)) \subseteq {\rm J}(A)$, $\nu^\ast(
{\rm S}(A/I)) \subseteq {\rm S}(A)$ and $\nu^\ast({\rm Z}(A/I)) = 
{\rm Z}(A) \cap Az \subseteq {\rm Z}(A)$; here ${\rm J}(A)$ denotes 
the (Jacobson) radical of $A$, and ${\rm S}(A)$ denotes the (left and
right) socle of $A$; see [1]. We also conclude:
\bigskip\noindent 
{\bf Proposition 2.} {\it In the situation above, let $L$ be an ideal
of ${\rm Z}(A/I)$. Then $\nu^\ast(L)$ is an ideal of ${\rm Z}(A)$.}
\bigskip\noindent 
{\it Proof.} This follows since $L$ is a ${\rm Z}(A)$-module and since
$\nu^\ast$ is a monomorphism of ${\rm Z}(A)$-modules. 
\bigskip\noindent 
In particular, we have $\nu^\ast({\rm R}(A/I)) \subseteq {\rm R}(A)$
where ${\rm R}(A) := {\rm Z}(A) \cap {\rm S}(A)$ denotes the Reynolds
ideal, an ideal in ${\rm Z}(A)$ (cf.\ [1] or [6]). 
Similarly, $\nu^\ast({\rm H}(A/I))$ and $\nu^\ast({\rm Z}_0(A/I))$
are ideals of ${\rm Z}(A)$; here ${\rm H}(A)$ denotes the Higman ideal,
an ideal in ${\rm Z}(A)$, and ${\rm Z}_0(A) = {\rm H}(A)^2$ is the 
sum of the blocks of $A$ which are simple $F$-algebras. Hence ${\rm 
Z}_0(A)$ is also an ideal of ${\rm Z}(A)$; see [2]. 

We obtain the following version of Kn\"orr's theorem in the context of
symmetric algebras:
\bigskip\noindent 
{\bf Corollary 3.} {\it In the situation above, the dimension of the 
ideal $\nu^\ast({\rm Z}_0(A/I))$ of ${\rm Z}(A)$ is the number of blocks
of $A/I$ which are simple $F$-algebras.} 
\bigskip\noindent
The inclusion ${\rm Z}_0(A) \subseteq {\rm H}(A) \subseteq {\rm R}(A)$
(cf.\ [2]) implies that, in Corollary 3, we have
$$\nu^\ast({\rm Z}_0(A/I)) \subseteq \nu^\ast({\rm R}(A/I)) \subseteq
{\rm R}(A).$$
\bigskip\bigskip\noindent
{\bf 2. Group algebras}
\bigskip\noindent 
In the following, let $F$ be an algebraically closed field of characteristic
$p>0$, and let $Q$ be a normal $p$-subgroup of a finite group $G$. We denote
the kernel of the canonical epimorphism $\nu: FG \longrightarrow F[G/Q]$ by
$I$. Then $I = FG \cdot {\rm J}(FQ) = {\rm J}(FQ) \cdot FG$, and we identify
$FG/I$ and $F[G/Q]$. 

The group algebra $FG$ is a symmetric $F$-algebra with symmetrizing linear
form 
$$\lambda: FG \longrightarrow F, \quad \sum_{g \in G} \alpha_g g \longmapsto
\alpha_1.$$
We denote the analogous symmetrizing linear form of $F[G/Q]$ by $\mu$. Then,
in the notation of Section 1, we have $z = Q^+$. 

It is well-known that ${\rm H}(FG) = {\rm Tr}_1^G(FG)$ is the ideal of 
$(FG)^G = {\rm Z}(FG)$ which, as a vector space over $F$, is spanned by the 
class sums of the conjugacy classes of $p$-defect zero in $G$. Furthermore,
we have
$${\rm Z}_0(FG) = {\rm H}(FG) \cdot G_p^+$$
(cf.\ [2]). We are now in a position to prove Theorem 1.
\bigskip\noindent 
{\it Proof (of Theorem 1).} Recall first that a block of $FG$ is a simple 
$F$-algebra if and only if it has defect zero. Hence, by Corollary 3, the 
number of blocks of defect zero in $F[G/Q]$ is the dimension of the ideal
$\nu^\ast({\rm Z}_0(F[G/Q])) = \nu^\ast({\rm H}(F[G/Q]) \cdot (G/Q)_p^+)$
of ${\rm Z}(FG)$. Let $T$ be a set of representatives (in $G$) for the 
cosets of $p$-elements in $G/Q$. Then $\nu(T^+) = (G/Q)_p^+$ and $Q^+T^+
= G_p^+$. Moreover, we have $\nu^\ast(gQ) = gQ^+$ and 
$$\nu^\ast(gQ \cdot (G/Q)_p^+) = \nu^\ast(gQ)T^+ = gQ^+T^+ = gG_p^+$$
for $g \in G$. Thus we obtain:
$$\eqalign{\nu^\ast({\rm H}(F[G/Q]) \cdot (G/Q)_p^+) 
&= \nu^\ast({\rm Tr}_1^{G/Q}(F[G/Q]) \cdot (G/Q)_p^+) \cr
&= \nu^\ast({\rm Tr}_1^{G/Q}(F[G/Q] \cdot (G/Q)_p^+)) \cr
&= {\rm Tr}_Q^G(\nu^\ast(F[G/Q] \cdot (G/Q)_p^+)) 
= {\rm Tr}_Q^G(FG \cdot G_p^+), \cr}$$
and Theorem 1 follows.
\bigskip\noindent
A short computation shows:
$$\eqalign{{\rm Tr}_Q^G(gG_p^+) 
&= {\rm Tr}_{{\rm C}_G(g)Q}^G({\rm Tr}_Q^{{\rm C}_G(g)Q}(gG_p^+)) \cr
&= {\rm Tr}_{{\rm C}_G(g)Q}^G(|{\rm C}_G(g)Q:Q| gG_p^+) \cr
&= |{\rm C}_G(g): {\rm C}_Q(g)| {\rm Tr}_{{\rm C}_G(g)Q}^G(gG_p^+)\cr}$$
for $g \in G$. Thus we have ${\rm Tr}_Q^G(gG_p^+) = 0$ unless 
${\rm C}_Q(g)$ is a Sylow $p$-subgroup of ${\rm C}_G(g)$.

Similarly, if $0 \ne {\rm Tr}_Q^G(g \cdot G_p^+) = \nu^\ast({\rm Tr}_1^{G/Q}
(gQ) \cdot (G/Q)_p^+)$ then we have $0 \ne {\rm Tr}_1^{G/Q}(gQ) \cdot (G/Q)_p^+$. 
Thus Lemma N in [5] implies that there exists a Sylow $p$-subgroup $P$ of $G$
such that $P \cap gPg^{-1} = Q$.

We denote by $g_p$ and $g_{p'}$ the $p$-factor and the $p'$-factor of $g$.
If ${\rm C}_Q(g)$ is a Sylow $p$-subgroup of ${\rm C}_G(g)$ then $g_p \in
{\rm C}_G(g)_p = {\rm C}_Q(g)$ and therefore $g_pG_p^+ = G_p^+$. 
Consequently, we have ${\rm Tr}_Q^G(gG_p^+) = {\rm Tr}_Q^G(g_{p'} \cdot 
G_p^+)$. This shows that the vector space ${\rm Tr}_Q^G(FG \cdot G_p^+)$
is spanned by the elements ${\rm Tr}_Q^G(gG_p^+)$ where $g$ ranges over
the $p'$-elements in $G$ such that ${\rm C}_Q(g)$ is a Sylow $p$-subgroup
of ${\rm C}_G(g)$ and such that there is a Sylow $p$-subgroup $P$ of $G$
such that $P \cap gPg^{-1} = Q$.

On the other hand, we have 
$${\rm Tr}_Q^G(FG \cdot G_p^+) \subseteq {\rm Tr}_Q^G((FG)^Q) =: (FG)_Q^G.$$
It is well-known that $(FG)_Q^G$ is an ideal in ${\rm Z}(FG)$ which, as a
vector space over $F$, is spanned by the class sums of the conjugacy classes
of $G$ whose defect groups are contained in $Q$.
\bigskip\noindent
{\bf Acknowledgement.} The author is grateful to G.~R.~Robinson for 
discussions which led to the present paper.
\bigskip\bigskip\noindent
{\bf References}
\bigskip\noindent
\item{1.} S.~Brenner and B.~K\"ulshammer, Ideals in the center of symmetric
algebras, {\it Int.\ Electron.\ J.\ Algebra} {\bf 34} (2023), 126-151
\smallskip\noindent
\item{2.} L.~H\'ethelyi, E.~Horv\'ath, B.~K\"ulshammer and J.~Murray,
Central ideals and Cartan invariants of symmetric algebras, {\it J.\ 
Algebra} {\bf 293} (2005), 243-260
\smallskip\noindent
\item{3.} R.~Kn\"orr, Counting blocks of defect zero, {\it Progr.\ Math.}
{\bf 95} (1991), 419-424
\smallskip\noindent 
\item{4.} B.~K\"ulshammer, Bemerkungen \"uber die Gruppenalgebra als 
symmetrische Algebra II, {\it J.\ Algebra} {\bf 75} (1982), 59-69
\smallskip\noindent 
\item{5.} B.~K\"ulshammer, Bemerkungen \"uber die Gruppenalgebra als 
symmetrische Algebra III, {\it J.\ Algebra} {\bf 88} (1984), 273-291
\smallskip\noindent 
\item{6.} B.~ K\"ulshammer, Group-theoretical descriptions of ring-theoretical
invariants of group algebras, {\it Progr.\ Math.} {\bf 95} (1991), 
425-442
\smallskip\noindent 
\item{7.} G.~R.~Robinson, The number of blocks with a given defect group,
{\it J.\ Algebra} {\bf 84} (1983), 556-566
\smallskip\noindent
\item{8.} Y.~Tsushima, On the annihilator ideals of the radical of a group
algebra, {\it Osaka J.\ Math.} {\bf 8} (1971), 91-97
\smallskip\noindent
\item{9.} L.~Wang and J.~Zhang, Orbital characters and their applications,
{\it J.\ Pure Appl.\ Algebra} {\bf 219} (2015), 121-141
\end